\documentclass[12pt]{amsart}

\usepackage{amsmath,amsthm,amssymb}
\usepackage[T1]{fontenc}
\usepackage{latexsym}
\usepackage{subfigure, epsfig, epsf}
\usepackage{color,fancyhdr,lastpage,graphicx}
\usepackage{xspace}
\usepackage{setspace}
\usepackage{bbm}

\newcommand{\NN}{\mathbb{N}}      
\newcommand{\RR}{\mathbb{R}}

\newcommand{\mcC}{\mathcal{C}} \newcommand{\mcN}{\mathcal{N}}  \newcommand{\mcA}{\mathcal{A}}

\newcommand{\whP}{\widehat{P}}

\newcommand{\EE}{\mathbb{E}} \newcommand{\PP}{\mathbb{P}}

\newcommand{\bm}{Brownian motion }

\newcommand{\sde}{stochastic differential equation }

\newcommand{\LEB}{\textsc{Leb}}

\newcommand{\st}{such that }

\newcommand{\varep}{\varepsilon}
\newcommand{\ep}{\epsilon}

\renewcommand{\leq}{\leqslant}
\renewcommand{\geq}{\geqslant}

\newcommand{\al}{\alpha}
\newcommand{\ra}{\rightarrow}
\newcommand{\ssk}{\smallskip}
\newcommand{\wrt}{with respect to }
\newcommand{\noi}{\noindent }

\newtheorem{thm}{\hspace{-0.15cm} {\sc Theorem} }

\newtheorem{lem}[thm]{\hspace{-0.18cm} {\sc Lemma} }

\newtheorem{prop}[thm]{\hspace{-0.18cm} {\sc Proposition}}

\numberwithin{equation}{section} 

\newenvironment{Dem}{%
    \begin{list}{\hspace{0.6cm}{\sc Proof --}}{%
        \setlength{\topsep}{0pt}%
        \setlength{\leftmargin}{0pt}%
        \setlength{\rightmargin}{0pt}%
        \setlength{\listparindent}{0pt}%
        \setlength{\itemindent}{0pt}%
        \setlength{\parsep}{0pt}%
        \addtolength{\leftmargin}{20pt}%
        \addtolength{\rightmargin}{0pt}%
    } \item }{\hfill{\space $\rhd$}\end{list}\smallskip}

\newenvironment{DemPropLip}{%
    \begin{list}{\hspace{0.68cm}{\sc Proof of proposition \ref{PropLipschitz} --}}{%
        \setlength{\topsep}{0pt}%
        \setlength{\leftmargin}{0pt}%
        \setlength{\rightmargin}{0pt}%
        \setlength{\listparindent}{0pt}%
        \setlength{\itemindent}{0pt}%
        \setlength{\parsep}{0pt}%
        \addtolength{\leftmargin}{20pt}%
        \addtolength{\rightmargin}{0pt}%
    } \item }{\hfill{\space $\rhd$}\end{list}\smallskip}

\newenvironment{DemPropAllTimes}{%
    \begin{list}{\hspace{0.68cm} {\sc Proof of proposition \ref{PropAllTimes} --}}{%
        \setlength{\topsep}{0pt}%
        \setlength{\leftmargin}{0pt}%
        \setlength{\rightmargin}{0pt}%
        \setlength{\listparindent}{0pt}%
        \setlength{\itemindent}{0pt}%
        \setlength{\parsep}{0pt}%
        \addtolength{\leftmargin}{20pt}%
        \addtolength{\rightmargin}{0pt}%
    } \item }{\hfill{\space $\rhd$}\end{list}\smallskip}

\newenvironment{SousDem}{%
    \begin{list}{\hspace{0.2cm} {\sc Proof --}}{%
        \setlength{\topsep}{0pt}%
        \setlength{\leftmargin}{0pt}%
        \setlength{\rightmargin}{0pt}%
        \setlength{\listparindent}{0pt}%
        \setlength{\itemindent}{0pt}%
        \setlength{\parsep}{0pt}%
        \addtolength{\leftmargin}{15pt}%
        \addtolength{\rightmargin}{0pt}%
    } \item }{\hfill{\space $\square$}\end{list}\smallskip}

\title{{\sc Spatial coagulation with bounded coagulation rate}}
\date{\today}
\author[{\sc I.F. Bailleul}]{{\sc I.F. Bailleul}}
\address{Statistical Laboratory, Center for Mathematical Sciences, Wilberforce Road, Cambridge, CB3 0WB, UK}
\email{i.bailleul@statslab.cam.ac.uk}
\urladdr{http://www.statslab.cam.ac.uk/~ismael/}


\begin{document}

\maketitle

\begin{abstract}
We prove that the spatial coagulation equation with bounded coagulation rate is well-posed for all times in a given class of kernels if the convection term of the underlying particle dynamics has divergence bounded below by a positive constant. Multiple coagulations, fragmentation and scattering are also considered.
\end{abstract}

\section{Introduction}

This article is concerned with Smoluchowski's picture of a chemical reaction involving two kinds of dynamics. Particles involved in the reaction are described by their location $x\in\RR^n$ and an inner structure, which we describe as an element $y$ of an Abelian semi-group. It can be its mass, in $\RR^*$ or $\NN$, or a finer structure saying for instance which smaller masses compose that mass, as in \cite{James2}. On the one hand, particles move as a result of the action of some convection flow and thermal diffusive effects. On the other hand, close particles may coagulate to form a new particle with a new characteristic. We shall mainly be interested here in binary collisions. Not to overload notations, let us describe a particle by its position $x\in\RR^n$ and its mass $y\in\RR^*_+$.

\ssk

The evolution of the system is described by a time dependent kernel $\mu_t(x\,;dy)$ from $\RR^n$ to $\RR^*_+$, with $\mu_t\bigl(x\,;(\al,\beta)\bigr)$ representing the concentration of particles at location $x$ with a mass in between $\al$ and $\beta$. Let $K$ be a measurable symmetric function on $\RR^*_+\times\RR^*_+$. Given a signed measure $\mu$ on $\RR_+^*$, define a signed measure $K(\mu,\mu)$ on $\RR^*_+$ setting 
$$
K(\mu,\mu) = \frac{1}{2}\int\bigl\{\delta_{y+y'}-\delta_y-\delta_{y'}\bigr\}K(y,y')\,\mu(dy)\mu(dy').
$$
Smoluchowski's picture of the microscopic dynamics leads heuristically to the formal evolution equation on kernels
\begin{equation}
\label{EquationDiff}
\partial_t\mu_t = L_x^*\mu_t + K(\mu_t,\mu_t),
\end{equation}
where $L_x = \frac{1}{2}\,a^{ij}(x,y)\,\partial_{x_i}\partial_{x_j} + b^i(x,y)\,\partial_{x_i}$ is a convection-diffusion differential operator. This heuristic result was proved in an important special case by Norris in \cite{JamesBrownianCoagulation}, and later in \cite{HammondRezakhanlou} and \cite{YaghoutiRezakhanlouHammond} by Hammond, Rezakhanlou and Yaghouti. Well-posedness for equation \eqref{EquationDiff} has been investigated in a relatively small number of works. Amann's article \cite{Amann} is the only work where the problem is tackled in this generality. Supposing all the coefficients in the dynamics bounded, he obtains the existence of a unique maximal solution depending continuously on the initial condition using highly non-trivial functional analytic tools. Under some smoothness and moment hypothesis on the initial condition, he is able to prove well-posedness for all times in space dimension $1$ or if the diffusivity coefficients do not depend on $y$, and $b=0$. This work was partly extended in \cite{AmannWalker}, in the Brownian case where $b=0$ and 
\begin{equation}
\label{BrownianDiffusivity}
\frac{1}{2}\,a^{ij}(x,y)\,\partial_{x_i}\partial_{x_j} = d(y)\triangle_x,
\end{equation}
with $d(y)$ bounded away from $0$ and $\infty$, and well-posedness for all times was proved under some moment and size conditions on the initial condition. They allow linear growth for the coagulation rate. The works \cite{MischlerRodriguez}, \cite{LaurencotMischlerRevista}, \cite{LaurencotMischlerArchive}, \cite{CanizoDesvillettesFellner}, identify different settings in which one can prove the existence of \textit{a} weak solution to equation \eqref{EquationDiff} (in different appropriate senses), defined for all times and for unbounded coagulation kernels  controlled in different ways. They all investigate equation \eqref{EquationDiff} without convection and with a diffusion of the form \eqref{BrownianDiffusivity}. Hammond and Rezakhanlou are able to give in \cite{HammondRezakhanlouMoment} a well-posedness result for all times for a discrete mass space, a bounded decreasing diffusivity $d(y)$, and some moment conditions. (See also \cite{RezakhanlouProc}.) 

\ssk

All these works investigate the differential problem \eqref{EquationDiff} in different functional settings. Norris \cite{JamesBrownianCoagulation} and Laurencot-Mischler \cite{LaurencotMischlerRevista} consider a mild integral formulation of equation \eqref{EquationDiff} rather than its differential counterpart. They work with no convection term and a diffusivity of the form \eqref{BrownianDiffusivity}. Relying crucially on Gaussian estimates, Norris is able to prove a well-posedness result for all times under some weak control on the coagulation rate and some moment condition on the initial condition. His method does not apply to position dependent diffusivity dynamics. Consult also \cite{RezakhanlouProc} for a recent existence and uniqueness result under some moment assumptions.

\ssk

We refer the reader to the above works for further and older references on the spatial coagulation equation \eqref{EquationDiff}. Note yet that the methods developped in \cite{AmannWalker} and \cite{LaurencotMischlerArchive}, \cite{LaurencotMischlerRevista}, should be sufficiently robust to deal with more general dynamics than \eqref{BrownianDiffusivity} on bounded smooth domains as they essentially rely on some compactness property of the operator.

\bigskip

Following Amann's view \cite{Amann}, we consider in this work a bounded coagulation kernel and a general convection-diffusion spatial dynamics. It is quite likely that under some reasonably general conditions on the convection-diffusion coefficients equation \eqref{EquationDiff} is well-posed for all times in a good functional framework. Despite all the subtle works done, this remains yet an open problem. We identify in this work a previously unnoticed case where this problem can be solved. Roughly speaking, a condition on the divergence of the convection term prevents the particles in the underlying microscopic dynamics from concentrating too much; this enables us to obtain some uniform bound in the spatial coordinate. Although far from opening a road to solving the general problem, the simplicity of our framework and argument contrasts with the sophistications developped in the above references.

\ssk

 We suppose that the coefficients $a^{ij}$ and $b^i$ are measurable functions of $(x,y)\in\RR^n\times\RR^*_+$, of class $\mcC^\infty$ as functions of $x\in\RR^n$, and with all their derivatives bounded, uniformly in $y\in\RR^*_+$; this allows linear growth of the drift. To give an integral form of \eqref{EquationDiff}, let us suppose here that one can write the matrix $a$ as $\sigma\sigma^*$ for some bounded matrix-valued measurable function $\sigma$ on $\RR^n\times\RR^*_+$, of class $\mcC^\infty$ as a function of $x$.  Denote by $\phi_t(x,y)$ the well-defined random flow on $\RR^n\times\RR^*_+$ of the \sde
\begin{equation}
\label{EquationDynamics}
\begin{split}
&dx_t = \sigma(x_t,y_t)\,dB_t+b(x_t,y_t)dt, \\
&dy_t = 0.
\end{split}
\end{equation}
For a kernel $\nu$ on $\RR^n\times\RR^*_+$, and $t\geq 0$, set
$$
\bigr(\whP_t\nu\bigl)(x\,;dy) = \EE\Bigl[e^{-\int_0^t (\textrm{div}\,b)(\phi_s(x,y),y)ds}\,\nu\bigl(\phi_t(x,y)\,;dy\bigr)\Bigr], 
$$
where $\textrm{div}\,b$ is the divergence of the vector field $b(\cdot,y)$ on $\RR^n$, for fixed $y$. The integral form of equation \eqref{EquationDiff} is 
\begin{equation}
\label{IntegralEquation}
\mu_t = \widehat{P}_t\mu_0 + \int_0^t\widehat{P}_{t-s}\bigl(K(\mu_s,\mu_s)\bigr)\,ds.
\end{equation}

\bigskip

\noi We now introduce the functional setting in which we are going to solve this equation. Let $\nu_0$ be a non-negative finite measure on $\RR_+^*$ \st 
\begin{equation}
\label{ProprieteNu0}
\nu_0\star\nu_0\leq C\,\nu_0
\end{equation}
for some positive constant $C$. Define $\mcN$ as the set of kernels $\nu(x\,;dy)$ from $\RR^d$ to $\RR_+^*$ of the form $\bigl\{f(x;y)\,\nu_0(dy)\bigr\}_{x\in\RR^d}$, where $f(x\,;\,y)$ is a bounded measurable function \textit{uniformly equicontinuous in $x$}: 
$$
\forall\,\ep>0,\,\exists\,\delta>0,\,\forall\,y\in\RR^*_+,\;|x-x'|\leq\delta\Rightarrow \bigl|f(x\,;\,y)-f(x'\,;\,y)\bigr|\leq\ep.
$$
Setting 
$$
\|\nu\| = \sup_{x\in\RR^d}\underset{y\in\RR^*_+}{\textrm{$\nu_0$-ess}\sup}\,\bigl|f(x\,;\,y)\bigr|
$$
defines a norm on $\mcN$. 

\begin{thm}
\label{MainTheorem}
Suppose the following two conditions hold.
\begin{itemize}
   \item[(i)] There exists some positive constant $\al,\beta$ \st $\al\,\emph{Id}\leq a(x,y)\leq \beta\,\emph{Id}$, for all $(x,y)\in\RR^n\times\RR^*_+$.
   \item[(ii)] There exists a positive constant $\varep$ \st the vector $b(\cdot,y)$ has a divergence bounded below by $\varep$, for all $y\in\RR^*_+$. 
\end{itemize}
\noi Let $M>0$ be an upper bound for $K$. Then, for any $\mu_0\in\mcN$, equation \eqref{IntegralEquation} has a unique maximal solution in $\mcN$, continuous \wrt $t$, and defined for all times provided $\|\mu_0\|$ is small enough; it is non-negative if $\mu_0$ is non-negative.
\end{thm}

\ssk

\noi \textsc{Remarks.}
\begin{enumerate}
   \item Use Laplace transform to construct a positive continuous and  $\LEB$-integrable function $g$ on $\RR_+^*$ \st the finite measure $\nu_0(dy) := g(y)\,dy$ satisfies \eqref{ProprieteNu0} for some positive constant $C$. Suppose $\mu_0(x\,;dy) = h_0(x\,;y)\,dy$, for a bounded measurable function $h_0$ with support in $\RR^n\times (0,m)$, for some $m<\infty$. It follows from theorem \ref{MainTheorem} that $\mu_t$ is well-defined for all times provided $\|h_0\|_\infty$ is small enough, and each $\mu_t(x\,;dy)$ is absolutely continuous \wrt Lebesgue measure.

   \item It will be clear from the proof of theorem \ref{MainTheorem} that the conclusion remains unchanged if the coagulation rate $K(y,y')$ is allowed to depend on $x$ also, as an equicontinuous function with a modulus of continuity independent of $y, y'\in\RR^*_+$. 

   \item Although all our statements are given here and below with $y\in\RR^*_+$, all the results hold true with $y$ taking values in a general Abelian semi-group describing the inner structure of particles, at the price of introducing the formalism used in \cite{James2}. All results hold in particular for the discrete masses spatial coagulation equation.

   \item Note that since we are working with a bounded coagulation kernel we do not require any moment assumption on the initial condition.

   \item A similar result holds when $\RR^n$ is replaced by a bounded domain of $\RR^n$, with smooth boundary, and Dirichlet or Neumann conditions are imposed on the boundary. Condition (ii) makes unnecessary the use of spectral properties of the diffusion/convection operators.
 
   \item Our method of proof is similar in spirit to the method used in \cite{AmannWalker}, where the quadratic coagulation term is balanced by a linear term\footnote{The main ingredients of \cite{AmannWalker} are mass conservation of solutions and the exponential decay of a semigroup on functions with zero spatial mean. We do not need mass conservation but use crucially the exponential decay of $\widehat{P}_t$.}. Our functional setting is more elementary, and the above result new.

   \item The introduction of the map $\widehat{P}_t$ is better understood by considering first a pure convection dynamics, where $a$ (and so $\sigma$) is null. In that case, the equation $\dot\mu_t=L_x^*\mu_t$ on kernels reads 
\begin{equation}
\label{OdeCase}
\dot\mu_t = -b^i\partial_{x_i}\mu_t - (\textrm{div}\,b)\mu_t.
\end{equation}
Introducing the flow $\varphi_t$ of the ordinary differential equation $\dot x = b(x)$ on $\RR^n$, equation \eqref{OdeCase}  has solution 
$$
\mu_t(x\,;dy) = \mu_0\bigl(\varphi_t(x)\,;dy\bigr)\,e^{-\int_0^t(\textrm{div}\,b)(\varphi_s(x),y)ds}
$$
In the case of a general second order differential operator $L_x$ as above, the associated dynamics is no longer deterministic, but given by the stochastic differential equation \eqref{EquationDynamics}. Equation \eqref{IntegralEquation} is the integral form of \eqref{EquationDiff} obtained by varying the constant.
\end{enumerate}

\section{Proof of theorem \ref{MainTheorem}}

\subsection{Well-posedness}

The proof of the following lemma is straightforward and left to the reader.

\begin{lem}
The space $\bigl(\mcN,\|\cdot\|\bigr)$ is a Banach space. 
\end{lem}

The proof of theorem \ref{MainTheorem} is an application of the elementary fixed point theorem for contractions of a Banach space. Fix $T>0$. Given an $\mcN$-valued path $\nu=(\nu_s)_{0\leq s\leq T}$, define for each $t\in [0,T]$ the kernel
$$
F\bigl(\nu\bigr)_t = \whP_t\nu_0 + \int_0^t\whP_{t-s}\bigl(K(\nu_s, \nu_s)\bigr)ds.
$$

\begin{prop}
\label{PropLipschitz}
The map $F$ is a locally Lipschitz map from $\mcC\bigl([0,T],\mcN\bigr)$ into itself.
\end{prop}

\noi It follows classically from proposition \ref{PropLipschitz} that one can associate to any $\mu_0\in\mcN$ a time $T$ depending only on $\|\mu_0\|$, in an increasing way, \st $F$ has a unique fixed point in $\mcC\bigl([0,T],\mcN\bigr)$ with initial value $\mu_0$. The following proposition proves then that $T=\infty$.

\begin{prop}
\label{PropAllTimes}
We have $\|\mu_t\|\leq \|\mu_0\|$, for all $t\in [0,T]$, if $\|\mu_0\|$ is small enough.
\end{prop}

We now proceed to proving the two propositions. 

\bigskip

\begin{DemPropLip}
We first deal with the operator $\widehat{P}_t$, for $t>0$.
\begin{lem}
\label{LemmaWidehatP}
\noi \emph{\textbf{a)}} The linear map $\whP_t$ sends $\mcN$ into itself and has norm no greater than $e^{-\varep t}$.

\noi \emph{\textbf{b)}} Fix $\nu\in\mcN$. The map $t\in\RR_+\mapsto \whP_t\nu\in\bigl(\mcN,\|\cdot\|\bigr)$ is continuous.
\end{lem}

\begin{SousDem}
\noi \textbf{a)} For $\mu(x\,;dy) = f(x;y)\nu_0(dy)$ in $\mcN$, 
$$
\bigl(\widehat{P}_t\mu\bigr)(x\,;dy) = \EE\Bigl[e^{-\int_0^t(\textrm{div}\,b)(\phi_s(x,y),y)\,ds}f\bigl(\phi_t(x,y);y\bigr)\Bigr]\,\nu_0(dy),
$$
so we need to see that the above expectation defines a ($y$-)uniformly equi-continuous function of $x$; denote it by $\bigl(P_t^{\,\textrm{div}b}f\bigr)(x,y)$. Rewrite equation \eqref{EquationDynamics} as a \sde driven by vector fields $V_i$
$$
dx_t = V_i(x_t,y)\,dB^i_t + b(x_t)\,dt,
$$
with $V_i(x,y)$ being the $i^{\textrm{th}}$ column of the matrix $\sigma(x,y)$. Kusuoka and Stroock have proved in \cite{KusuokaStroockIII} that we have\footnote{They prove their results for $\sigma$ and $b$ and all their derivatives bounded; their proof of the following result extends in a straightforward way to our elliptic setting.} for each $i\in\{1,\dots,n\}$
\begin{equation}
\label{ResultKusuoka}
\underset{x\in\RR^n}{\sup}\,\bigl|V_iP_t^{\,\textrm{div}b}f\bigr|(x,y) \leq c\,t^{-\frac{1}{2}}\,\underset{x\in\RR^n}{\sup}\,\bigl|f(x\,;\,y)\bigr|,
\end{equation}
for each $y\in\RR^*_+$, where the constant $c$ depends only on an upper bound for the uniform norm of the $V_i$'s and all their derivatives. It can be taken to be independent of $y\in\RR^*_+$ from our assumptions. Also, from the uniform ellipticity of the matrices $a(x,y)$, inequality \eqref{ResultKusuoka} implies\footnote{Due to the uniform ellipticity of the matrices $a(x,y)$, inequality \eqref{ResultKusuoka} is equivalent to the existence of a constant $c'$, independent of $i\in\{1,\dots,n\}$, \st $\underset{x\in\RR^n}{\sup}\,\bigl|\partial_{x_i}P_t^{\textrm{div}b}f\bigr|\leq c't^{-1/2}\underset{x\in\RR^n}{\sup}\,\bigl|f(x,y)\bigr|$. The $y$-uniform equi-continuity of $P_t^{\textrm{div}b}f$ follows from the boundedness of $f$ as $t>0$ is fixed.} the ($y$-)uniform equicontinuity of $P_t^{\,\textrm{div}b}f$. This proves that $\widehat{P}_t$ sends $\mcN$ into itself; the statement on its norm is straightforward from hypothesis (ii) of theorem \ref{MainTheorem}.

\ssk

\noi \textbf{b)} It suffices to prove the continuity at $t=0$. For $\nu(x,dy) = f(x\,;\,y)\nu_0(dy)$ in $\mcN$, we have
\begin{equation}
\label{IneqContinuitywhatPt}
\begin{split}
\|\whP_t\nu-\nu\| &\leq \sup_{x,y}\,\EE\Bigl[\Bigl|e^{-\int_0^t(\textrm{div}\,b)(\phi_s(x,y),y)\,ds}f(\phi_t(x,y)\,;y)-f(x\,;\,y)\Bigr|\Bigr] \\
                  &\leq \sup_{x,y}\,\EE\Bigl[\bigl|f(\phi_t(x,y)\,;y)-f(x\,;\,y)\bigr|\Bigr] + \|f\|_\infty\,\sup_{x,y}\,\EE\Bigl[\Bigl|e^{-\int_0^t (\textrm{div}\,b)(\phi_s(x,y)\,;y)\,ds}-1\Bigr|\Bigr] \\
                  &\leq \sup_{x,y}\,\EE\Bigl[\bigl|f(\phi_t(x,y)\,;y)-f(x\,;\,y)\bigr|\Bigr] + c\,\|f\|_\infty t.
\end{split}
\end{equation}
for some positive constant $c$. By the uniform equi-continuity of $f(\cdot\,;y)$ one can associate to any $\al>0$ a $\delta>0$ \st
$$
\forall\,y\in\RR^*_+, \; \forall\,x,x'\in\RR^n, \;\;|x'-x|\leq \delta \Rightarrow \big|f(x,;y)-f(x'\,;y)\big|\leq \al.
$$
Then, we have for any $(x,y)\in\RR^n\times\RR^*_+$
\begin{equation*}
\begin{split}
\EE\bigl[\big|f(\phi_t(x,y)\,;y)-f(x\,;y)\big|\bigr] &\leq \al\,\PP\bigl(|\phi_t(x,y)-x|\leq \delta\bigr) + 2\|f\|_\infty\,\PP\bigl(|\phi_t(x,y)-x|>\delta\bigr) \\
&\leq \al + 2\|f\|_\infty\,\frac{\EE\bigl[|\phi_t(x,y)-x|^2\bigr]}{\delta^2}.
\end{split}
\end{equation*}
As the isometry property of stochastic integration \wrt \bm and Jensen's inequality give
\begin{equation*}
\begin{split}
\EE\bigl[|\phi_t(x,y)-x|^2\bigr] & \leq 2\EE\Bigl[\Bigl(\int_0^t a(\phi_s(x,y),y)dB_s\Bigr)^2\Bigr] + 2\EE\Bigl[\Bigl(\int_0^t b(\phi_s(x,y),y)ds\Bigr)^2\Bigr] \\ 
& \leq 2\,\EE\Bigl[\int_0^t\bigl|a(\phi_s(x,y),y)\bigr|^2\,ds\Bigr] + 2t\,\EE\Bigl[\int_0^t\bigl|b(\phi_s(x,y),y)\bigr|^2\,ds\Bigr] \leq Ct
\end{split}
\end{equation*}
for some positive constant $C$ depending only on $\|a\|_\infty$ and $\|b\|_\infty$, it follows that 
$$
\sup_{(x,y)\in\RR^n\times\RR^*_+}\,\EE\Bigl[\bigl|f(\phi_t(x,y)\,;y)-f(x\,;\,y)\bigr|\Bigr] \leq \al+\frac{2C\|f\|_\infty}{\delta^2}\,t.
$$
\noi Together with \eqref{IneqContinuitywhatPt}, this inequality implies that $\|\whP_t\nu-\nu\| \ra 0$ as $t$ goes to $0$.
\end{SousDem}

\medskip

It follows from lemma \ref{LemmaWidehatP} that proposition \ref{PropLipschitz} will be proved if we can prove that the kernel $K(\mu,\mu)$ belongs $\mcN$ if $\mu$ does (\textbf{c}), and is a locally Lipschitz function of $\mu$ (\textbf{d}). Recall we write $M$ for an upper bound of $K$.

\medskip

\noi \textbf{c)} Given $\mu(x,dy) = f(x\,;\,y)\,\nu_0(dy)$ in $\mcN$, set, for any $x\in\RR^n$,
$$
K^+(\mu,\mu)(x\,;dz) = \frac{1}{2}\int \delta_{y+y'}(dz)\,K(y,y')\,f(x\,;y')\,\nu_0(dy')\,f(x\,;\,y)\nu_0(dy),
$$
and 
$$
K^-(\mu,\mu)(x\,;dz) = \Bigl(\int K(y,z)f(x\,;y)\,\nu_0(dy)\Bigr)\,f(x;z)\,\nu_0(dz),
$$
so that $K(\mu,\mu) = K^+(\mu,\mu)-K^-(\mu,\mu)$. As $\nu_0$ satisfies hypothesis \eqref{ProprieteNu0} we have\footnote{We write $|\rho|$ for the absolute value of a signed measure $\rho$.}
$$
\bigl|K^+(\mu,\mu)(x\,;\cdot)\bigr| \leq M\|f\|_{\infty}^2\,\nu_0\star\nu_0 \leq  CM\|f\|_{\infty}^2\,\nu_0;
$$
we can thus write
$$
K^+(\mu,\mu)(x\,;dy) = k^+(x,y)\,\nu_0(dy)
$$
for some bounded measurable function $k^+(x,y)$. Taking any bounded function $g$, we see by dominated convergence that the integral 
$$
\int g(y)\,K^+(\mu,\mu)(x\,;dy)
$$ 
is a uniformly continuous function of $x$, with a modulus of continuity independent of $g$. It follows that one can construct a version of $k^+(x,y)$ which is ($y$-)uniformly equicontinuous in $x$. So $K^+$ sends $\mcN$ into itself; it is straightforward to see that $K^-$ does the same.

\ssk

\noi \textbf{d)} Denote by $|\cdot|_{{\sc TV}}$ the total variation norm on finite signed measures on $\RR^*_+$. The elementary inequality  
$$
\bigl|\bigl(K(\mu,\mu)-K(\mu',\mu')\bigr)(x\,;\cdot)\bigr|_{\sc{TV}} \leq M\Bigl(|\mu(x\,;\cdot)|_{\sc{TV}}+|\mu'(x\,;\cdot)|_{\sc{TV}}\Bigr)\,\big|\,\mu(x\,;\cdot)-\mu'(x\,;\cdot)\,\big|_{\sc{TV}}, \quad x\in\RR^d,
$$
shows that $K$ is locally Lipschitz in $\bigl(\mcN,\|\cdot\|\bigr)$.
\end{DemPropLip}

\medskip

\begin{DemPropAllTimes}
Set $m = \nu_0(\RR^*_+)$. As equation \eqref{IntegralEquation} gives for all $t\in [0,T]$
$$
\|\mu_t\| \leq e^{-\varep t}\|\mu_0\|+M\Bigl(\frac{C}{2}+m\Bigr)\int_0^t e^{-\varep (t-s)}\|\mu_s\|^2\,ds,
$$
the real-valued function $\|\mu_t\|$ is bounded above by the solution of the equation 
$$
\dot{z}_t = -\varep z_t + M\Bigl(\frac{C}{2}+m\Bigr)\,z_t^2, \quad z_0=\|\mu_0\|, \quad 0\leq t\leq T.
$$
Its solution 
$$
z_t = \frac{\varep}{M\bigl(\frac{C}{2}+m\bigr)+\frac{\varep-M\bigl(\frac{C}{2}+m\bigr)\|\mu_0\|}{\|\mu_0\|}e^{\varep\,t}}
$$
is defined for all times and  less than or equal to $z_0$ for $\|\mu_0\|<\frac{\varep}{M\bigl(\frac{C}{2}+m\bigr)}$.
\end{DemPropAllTimes}

\subsection{Positivity}

We do not suppose in this section that $\mu_0\in\mcN$ is small, so we work on a time interval $[0,T]$ where equation \eqref{IntegralEquation} is known to have a unique solution in $\mcN$ continuous \wrt $t$. 

\ssk

Suppose $\mu_0\geq 0$. We prove the positivity of $\mu_t$ by the usual method; see for instance theorem 10 in \cite{AmannWalker}. We recall it here for the reader's convenience. Let $\al$ be a positive constant \st 
\begin{equation}
\label{ConditionPositivity}
\int_{\RR^*_+}\mu_s(x\,;dy)\leq \|\mu_s\|\,\nu_0(\RR^*_+) \leq \frac{\al}{M}
\end{equation}
for all $x\in\RR^*_+$ and $s\in [0,T]$. Define a map $G$ from $\mcC\bigl([0,T],\mcN\bigr)$ into itself setting 
$$
G\bigl((\nu_s)_{0\leq s\leq T}\bigr)_t = \widehat{P}_t\mu_0 + \int_0^t\widehat{P}_{t-s}\bigl(K(\mu_s,\nu_s) + \al\,\nu_s\bigr)\,ds - \al\int_0^t\widehat{P}_{t-s}\nu_s\,ds.
$$
Given two measures $\rho$ and $\rho'$ on $\RR^*_+$, set 
$$
K(\rho,\rho')(dz) = \int\bigl\{\delta_{y+y'}(dz)-\delta_y(dz)-\delta_{y'}(dz)\bigr\}K(y,y')\,\rho(dy)\rho'(dy').
$$ 
As the following inequality between kernels holds pointwise in $x\in\RR^n$
\begin{equation*}
\begin{split}
K(\mu_s,\nu_s)(dz)&\geq -\Bigl(\int_{\RR^*_+}K(y,z)\,\mu_s(dy)\Bigr)\nu_s(dz) + \al\,\nu_s(dz) \\
                             &\geq \Bigl(\al - \int_{\RR^*_+}K(y,z)\,\mu_s(dy)\Bigr)\nu_s(dz)\overset{\eqref{ConditionPositivity}}{\geq} 0,
\end{split}
\end{equation*}
the first integral in the definition of $G$ is non-negative. The path $(\mu_s)_{0\leq s\leq T}$ solves by construction the equation on kernels
$$
\nu_t + \al\int_0^t\widehat{P}_{t-s}\nu_s\,ds = G\bigl((\nu_s)_{0\leq s\leq T}\bigr)_t \;\;\geq 0, \quad t\in [0,T].
$$
Taking $T$ smaller if necessary, one can solve this equation by a fixed point method; see e.g. theorem 10 of \cite{AmannWalker}. As each iteration step of the method preserves non-negative kernels, its unique fixed point $(\mu_s)_{0\leq s\leq T}$ is non-negative if $\mu_0$ is non-negative.

\bigskip

\noi \textsc{Remark on the proof of theorem \ref{MainTheorem}.} The pattern of proof used to prove theorem \ref{MainTheorem} works equally well for some sub-elliptic convection-diffusion operators $L_x$ of H\"ormander form $\frac{1}{2}\sum_{i=1}^pV_i^2+B$ satisfying the following assumptions.
\begin{itemize}
   \item The vector fields $V_i(\cdot\,;y)$ and $B(\cdot\,;y)$ on $\RR^n$ are measurable functions of $(x,y)\in\RR^n\times\RR_+^*$, of class $\mcC^\infty$ as functions of $x\in\RR^n$, with all their derivatives bounded uniformly in $y\in\RR^*_+$,
   \item The vector fields $B(\cdot,y)$ satisfy hypothesis (ii) of theorem \ref{MainTheorem},
   \item Set $\mcA = \{\emptyset\}\cup\bigcup_{k=1}^\infty\{1,\dots,p\}^k$, set $V_{[\emptyset]}=0$ and $V_{[i]}=V_i$, for $i=1\dots,p$, and define inductively the vector fields $V_{[\al]}$, for $\al\in\mcA$, setting
\begin{equation*}
V_{[\al i]}=\bigl[V_{[\al]},V_i\bigr], \;\textrm{ for }i=1,\dots,p.
\end{equation*}
Given a positive integer $\ell$, write $\mcA_\ell$ for $\{\emptyset\}\cup\bigcup_{k=1}^\ell\{1,\dots,p\}^k$. We suppose there exists an $\ell$ \st 
\begin{equation}
\label{UH}
\underset{\xi\in\mathbb{S}^{n-1}}{\inf}\,\Bigl\{\sum_{\al\in\mcA_\ell}\bigl(V_{[\al](x,y)},\xi\bigr)^2\,;\,(x,y)\in\RR^n\times\RR^*_+\Bigr\}>0.
\end{equation}
\end{itemize}
Indeed, theorem 2 of \cite{KusuokaMalliavinRevisited} provides a uniform in $\al\in\mcA_\ell$ estimate of 
$$
\underset{(x,y)\in\RR^n\times\RR^*_+}{\sup}\,\bigl|V_{[\al]}P_t^{\textrm{ div}b}f\bigr|(x,y)
$$ 
in terms of $\|f\|_\infty$, at a fixed time $t>0$, similar to \eqref{ResultKusuoka}. The $y$-uniform equi-continuity of $P_t^{\textrm{ div}\,b}f$ follows from condition \eqref{UH} as above. The remainder of the proof of theorem \ref{MainTheorem} does not need any change.

\ssk

This kind of improvement of theorem \ref{MainTheorem} seems difficult to get by the analytical methods developped in the works mentionned in the introduction.

\section{Multiple coagulations, fragmentation and scattering}

\subsection{Multiple coagulations}

It is clear from the proof of theorem \ref{MainTheorem} that it can be extended to a multiple coagulation framework. Write $K_2$ for the above binary coagulation kernel $K$ and let $K_n$ be for all $n\geq 3$ a symmetric non-negative valued function on $\bigl(\RR^*_+\bigr)^n$. Set 
$$
K_n(\mu,\cdots,\mu) = \frac{1}{n!} \int\bigl\{\delta_{y_1+\cdots+y_n}-(\delta_{y_1}+\cdots+\delta_{y_n})\bigr\}K_n(y_1,\cdots,y_n)\,\mu(dy_1)\cdots\mu(dy_n)
$$
The spatial coagulation equation takes in that general framework the form
\begin{equation}
\label{EquationMultipleCoagulation}
\partial_t\mu_t = L^*_x\mu_t + \sum_{n\geq 2}K_n(\mu_t,\dots,\mu_t).
\end{equation}
 
\begin{thm}
Suppose the coagulation kernels $K_n$ are uniformly bounded above by a constant independent of $n$. Suppose hypotheses (i) and (ii) of theorem \ref{MainTheorem} hold. Then equation \eqref{EquationMultipleCoagulation} has a unique maximal solution in $\mcN$, defined for all times if $\|\mu_0\|$ is small enough; it is non-negative if $\mu_0$ is non-negative.
\end{thm}

\begin{Dem}
The proof of this statement is similar to the proof of theorem \ref{MainTheorem}, noticing that we have $\nu_0^{*n}\leq C^n\,\nu_0$ under condition \eqref{ProprieteNu0}, and introducing, for any measure $\mu$ on $\RR^*_+$, the measures $K_n^+(\mu,\dots,\mu)$ and $K_n^-(\mu,\dots,\mu)$ defined by the formulas
\begin{equation*}
\begin{split}
&K_n^+(\mu,\dots,\mu)(dz) = \frac{1}{n!} \int\delta_{y_1+\cdots+y_n}(dz)K_n(y_1,\cdots,y_n)\,\mu(dy_1)\cdots\mu(dy_n), \\
&K_n^-(\mu,\dots,\mu)(dz) = \frac{1}{(n-1)!}\,\Bigl(\int K(y_1,\dots,y_{n-1},z)\,\mu(dy_1)\cdots\mu(dy_{n-1})\Bigr)\mu(dz).
\end{split}
\end{equation*} 
Denote by $M$ an upper bound for all the coagulation kernels, and write $m$ for $\nu_0(\RR^*_+)$. To prove that the unique maximal solution is defined for all times for a small enough initial condition, remark that equation \eqref{EquationMultipleCoagulation} implies the inequality
$$
\|\mu_t\| \leq e^{-\varep t}\|\mu_0\| + \sum_{n\geq 2}\Bigl(\frac{MC^n}{n!} +\frac{m^{n-1}}{(n-1)!}\Bigr)\int_0^te^{-\varep (t-s)}\|\mu_s\|^n\,ds,
$$
so $\|\mu_t\|$ is bounded above by the solution of the differential equation
$$
\dot{z}_t = -\varep z_t + M(e^{Cz_t}-1-Cz_t)+ z_t\bigl(e^{mz_t}-1\bigr), \quad z_0 = \|\mu_0\|.
$$
Its solution is decreasing and defined for all times provided $z_0$ is small enough.
\end{Dem}

\subsection{Coagulation, fragmentation and scattering}

More realistic models of chemical reactions have fragmentation and scattering terms. The fragmentation term accounts for the spontaneous breakage of particles into smaller particles; it is defined by the following formula, where $\rho$, here and below, is used as a generic measure on $\RR^*_+$,
$$
\frak{B}(\rho)(dz):= \int\bigl\{F(y\,;dz)-\delta_y(dz)\bigr\}B(y)\,\rho(dy).
$$
A particle of mass $y$ splits into smaller particles at rate $B(y)$, the result of the splitting being given by the non-negative measure $F(y\,;dz)$ with support in $(0,y)$ and total mass $\int z\,F(y\,;dz)=y$. The scattering term accounts for the fact that real particles have a mass no bigger than some maximum mass $y_0$, in a given situation, so if coagulation forms a particle with mass bigger than $y_0$, this particle breaks instantaneously into smaller particles. The scattering term $\frak{S}$ is determined by a scattering kernel $S$ from $(y_0,2y_0)$ to the set of non-negative measures on $\RR^*_+$, and the formula
$$
\frak{S}(\rho,\rho)(dz):= \int_{(0,y_0)^2}\bigl\{S(y+y';dz)-\delta_y(dz)-\delta_{y'}(dz)\bigr\}{\bf 1}_{y+y'>y_0}K(y,y')\,\rho(dy)\rho(dy').
$$
This term was first introduced in \cite{FassanoRosso}, and studied further in \cite{WalkerBreakage} and \cite{WalkerBreakageDiffusion}. The coagulation term is changed in this model into
$$
K(\rho,\rho) = \int\bigl\{\delta_{y+y'}-\delta_y-\delta_{y'}\bigr\}{\bf 1}_{y+y'\leq y_0}\,K(y,y')\,\rho(dy)\,\rho(dy').
$$
The proof of the following theorem is analogue to the proof of theorem \ref{MainTheorem}.

\begin{thm}
Let $y_0>0$ be an upper bound for the size of the particles in the system. Suppose the following conditions hold.
\begin{itemize}
   \item[(i)] There exists some positive constant $\al,\beta$ \st $\al\emph{Id}\leq a(x,y)\leq \beta\emph{Id}$, for all $(x,y)\in\RR^n\times (0,y_0)$.
   \item[(ii)] There exists a positive constant $\varep$ \st the vector field $b(\cdot,y)$ has divergence bounded below by $\varep$, for all $y\in (0,y_0)$.
   \item[(iii)] The coagulation and breakage rates $K, B$ are bounded, and the fragmentation kernel $F(y\,;dz)$ has the form $\frak{f}(y\,;z)\,\nu_0(dz)$, with $\|B\|_\infty\bigl(1+\|\frak{f}\|_\infty\bigr)<\varep$.
   \item[(iv)] We have 
\begin{equation}
\label{ConditionScattering}
\int S(y+y'\,;\,dz)\,\nu_0(dy)\nu_0(dy')\leq C\,\nu_0(dz)
\end{equation}
for some positive constant $C$.
\end{itemize}
Given $\mu_0\in\mcN$, the equation 
$$
\mu_t = \whP_t\mu_0 + \int_0^t\whP_{t-s}\bigl(K(\mu_s,\mu_s) + \frak{B}(\mu_s) + \frak{S}(\mu_s,\mu_s)\bigr)ds, \quad t\in[0,T]
$$
has a unique solution in $\mcN$ continuous \wrt $t$, and defined for all times provided $\|\mu_0\|$ is small enough; it is non-negative if $\mu_0$ is non-negative.
\end{thm}

Condition (iv) will for instance hold if $\nu_0=\LEB_{|(0,y_0)}$ and the measures $S(a\,;\,dz)$ have a uniformly bounded density \wrt Lebesgue measure on $(0,y_0)$, for $a\in (y_0,2y_0)$.

\end{document}